\documentclass[11pt,a4paper]{article}
\usepackage{amsfonts,amsmath,amstext,amssymb}
\usepackage{dsfont}
\usepackage{theorem}    
\usepackage{a4wide}
\usepackage{verbatim}

\def\R{\mathbb R}     

\newtheorem{lema}{Lemma}

\newtheorem{teor}[lema]{\bf Theorem}
\newtheorem{coro}[lema]{\bf Corollary}

\newtheorem{rema}[lema]{\bf Remark}
\newtheorem{counter}[lema]{\bf Counterexample}

\theorembodyfont{\rmfamily}

\renewcommand{\inf}{\mathrm{Inf}}
\renewcommand{\sup}{\mathrm{Sup}}

\title{New examples of Moser-Bernstein type problems for\\ some nonlinear elliptic partial differential\\ equations arising in geometry}

\author{Alfonso Romero${}^{a}$\footnote{This work has been partially supported by the Spanish MIINECO aud ERDF project MTM2016-78807-C2-1-P. The first author has been also partially supported by the Andalusian and ERDF project A-FQM-494-UGR18.},
Rafael M. Rubio${}^{b*}$ and  Juan J. Salamanca${}^{b*}$\\[6mm]
${}^a$ Departamento de Geometr\'{\i}a y Topolog\'{\i}a, \\
Universidad de Granada, 18071 Granada, Spain \\ E-mail\textup{:
\texttt{aromero@ugr.es}}\\[3mm]
${}^b$ Departamento de Matem\'aticas, Campus de Rabanales, \\[0.5mm]
Universidad de C\'ordoba, 14071 C\'ordoba, Spain,\\[0.5mm]
E-mails\textup{:\texttt{\;rmrubio@uco.es},\,\,\texttt{jjsalamanca@uco.es}}}

\date{}

\begin{document}

\maketitle

\thispagestyle{empty}

\begin{abstract}
A family of nonlinear partial differential equations of divergence form is considered. Each one is the Euler-Lagrange equation of a natural Riemaniann variational problem of geometric interest. New uniqueness results for the entire solutions of these equations on a parabolic Riemaniann manifold of arbitrary dimension are given. In particular, several Moser-Bernstein type theorems are proved.
\end{abstract}

\vspace{1,5mm}


\hyphenation{ma-ni-fold}

\section{Introduction}
Among the elliptic quasi-linear PDEs, the equation of minimal hypersurfaces in Euclidean
space 
\begin{equation}\label{edp2}
\mathrm{div} \left(  \frac{D u}{ \sqrt{1 + |D u|^2}} \right) =0 ,
\end{equation}
has a long and fruitful history and has deserved the attention of many researchers. From a geometric viewpoint, it is the Euler-Lagrange equation of a classical variational problem. In fact, for each $u\in C^{\infty}(\Omega)$, $\Omega$ an open domain in $\R^n$, the $n$-form $\sqrt{1+ |D u|^2}\,dV$ on
$\Omega$ represents the volume element of the induced metric from $\R^{n+1}$ on the graph $\Sigma_u=\{(u(x),x)\, : \, x\in \Omega\}$ and the critical points of the $n$-volume functional $u \mapsto \int \sqrt{1+ |D u|^2}\,dV$ are characterized by equation (\ref{edp2}).

\vspace{1mm}

The early seminal result of S. Bernstein in 1914 for $n=2$, \cite{Bernstein}, amended by E. Hopf in 1950, \cite{Hopf}, is the well-known uniqueness theorem,

\begin{quote}{\it
The only entire solutions to the minimal surface equation in $\R^3$ are the affine functions $$u(x,y)=a \, x+ b\, y+ c,$$ where $a,b,c\in \R$.}
\end{quote}

Actually, Bernstein obtained his result as an application of the so called Bernstein's geometric theorem,

\begin{quote}
{\it If the Gauss curvature of the graph of $u\in C^{\infty}(\R^2)$ in $\R^3$ satisfies  $K\leq 0$ everywhere and $K < 0$ at some point, then $u$ cannot be bounded.}
\end{quote}

As an application, Bernstein proved a very general Liouville theorem,

\begin{quote}
{\it Any bounded solution $u\in C^{\infty}(\R^2)$ of the equation
$$A \, u_{xx} + 2\,B \, u_{xy} + C \, u_{yy} = 0,$$
\noindent where $A,B,C\in C^{\infty}(\R^ 2)$ such that $AC-B^2 > 0$, must be constant.}
\end{quote}

The possible extension of the classical Bernstein theorem to higher dimension is known as the \emph{Bernstein conjecture}. It has been an amusing research topic for a long time and it made many
advances on geometric analysis (see \cite{Osserman} for a detailed survey until 1984).

\vspace{1mm}

A remarkable contribution to the Bernstein conjecture in 1961, due to J. Moser \cite{Moser},
was the following general result,

\begin{quote}
{\it The only entire solutions $u$ to the minimal surface equation in $\R^{n+1}$ such that $\vert Du\vert \leq C$, for some $C\in \R^+$, are the affine functions $$u(x_1, \ldots ,x_n)=a_1 x_1+ \ldots +a_n x_n + c,$$ 
$a_i,c\in \R$, $1\leq i \leq n$, with $\displaystyle\sum_{i=1}^na_i^2\leq C^2$.}
\end{quote}

This theorem is called the Moser's weak Bernstein theorem, or the \emph{Moser-Bernstein theorem} in short. On the other hand, L. Bers proved in 1951, \cite{Bers}, that a solution $u$ of the minimal surface equation in $\R^3$ defined on the exterior of a closed disc in $\R^2$ has bounded $\vert Du \vert$. Therefore, the Moser-Bernstein theorem for $n=2$ and Bers' result provided another proof of the Bernstein theorem.

\vspace{1mm}

In 1968, J. Simons \cite{Simons} proved a result which in combination with theorems of E. De Giorgi \cite{Giorgi} and W.H. Fleming \cite{Fleming} yield a proof of the Bernstein conjecture
for $n\leq 7$. Moreover, there is a counterexample $u\in C^{\infty}(\R^n)$ to the Bernstein conjecture for each $n\geq 8$, (of course, with unbounded $\vert Du \vert$).

\vspace{1mm}

For a Riemannian product $\mathbb{R}\times M^n$, a graph $\{(u(x),x)\, :\, x\in  M^n\}$ is minimal if and only if the function $u$ satisfies formally the same differential equation (\ref{edp2}), i.e., 
\begin{equation}\label{edp2.1}
\mathrm{div} \left(  \frac{D u}{ \sqrt{1 + |D u|^2}} \right) =0 ,
\end{equation}
\noindent where $Du$ is the gradient of $u$, $\mid Du\mid$ its length and ${\rm div}$ the divergence in the Riemannian manifold $(M^n,g)$.

In the context of $3$-dimensional Riemannian product spaces $\mathbb{R}\times M^2$, where $M^2$ denotes a complete Riemannian surface with non-negative Gaussian curvature, H.  Rosenberg \cite{Rosenberg_2002} showed  that an entire minimal graph in $\mathbb{R}\times M^2$ must be totally geodesic. In \cite{Al-Da}, L.J. Al\'\i{}as, M. Dajczer and J. Ripoll completed Rosenberg's result showing that an entire minimal graph in $\mathbb{R}\times M^2$, with $M^2$ a complete Riemannian surface with non-negative Gaussian curvature $K$, and $K(p_0)>0$ at some point $p_0 \in M^2$, must be a slice $\{t_0\}\times M^2$, $t_0\in \mathbb{R}$.

On the other hand, H. Rosenberg, F. Schulze and J. Spruck \cite{Rosenberg_2013}, shown that an entire
minimal graph with non-negative height function in a Riemannian product $\mathbb{R}\times M^n$, such
that $M^n$ is complete with non-negative Ricci curvature and sectional curvature bounded
from below, must be a slice $\{t_0\}\times M^n$, $t_0\in \mathbb{R}$.

More recently, A.M.S. Oliveira and H.F. de Lima \cite{O-L} obtain a new Moser-Bernstein
type result for a Riemannian product $\mathbb{R}\times M^n $, where $M^n$ is complete with non-negative Ricci curvature and sectional curvature bounded from below, under the additional assump
tion of boundedness of the norm of the second fundamental form of the entire graph of a
function $u \in  C^\infty (M^n)$, whose gradient is bounded. So, they shown that if the graph has
constant mean curvature, then it is minimal. Moreover, if the function $u$ is bounded from below,
then the corresponding graph must be a slice $\{t_0\}\times M^n$, $t_0\in \mathbb{R}$.

The main aim of this paper is to prove several Moser-Bernstein type results in ambient
Riemannian manifolds more general than a product Riemannian space $\mathbb{R}\times M^n$. Namely,
we consider certain warped Riemannian spaces defined as follows \cite{ON}: given a positive smooth function $f$ on $I$, consider the \emph{warped product} with \emph{base} $(I,dt^2)$, {\it fiber} $(M,g)$ and \emph{warping function} $f$, i.e., the product manifold $I\times M$, endowed with the
Riemannian metric
\begin{equation}
\bar g=\pi_{_I}^*(dt^2)+ f(\pi_{_I})^2 \pi_{_M}^*(g),
\end{equation}
\noindent where $\pi_{_I}$ and $\pi_{_M}$ denote the projections onto $I$ and $M$, respectively. Following the terminology of \cite{ON}, let us denote this Riemannian manifold by $I \times_f M$.

For each $u\in C^{\infty}(M)$ such that $u(M)\subset I$, the graph of $u$ in the Riemannian warped product $I\times_f M$ defines a hypersurface which is minimal if and only if $u$ satisfies
\begin{equation}\label{edp}
\mathrm{div} \left(\frac{D u}{f(u) \sqrt{f(u)^2 + |D u|^2}}
\right) = \frac{f'(u)}{\sqrt{f(u)^2 + |D u|^2}} \left\{\, n \,-\,
\frac{|D u|^2}{f(u)^2} \,\right\} \, ,
\end{equation}
which is a non-linear elliptic equation of divergence form. In accordance to the classical terminology, we will refer it as the minimal surface (MS) equation in $I \times_f M$. Note that if
$M=\R^n$ and $f=1$, then the MS equation (\ref{edp}) agrees to the classical MS equation (\ref{edp2}). 

Here, we are mainly interested in the case that $I \times_f M$ is far from a Riemannian product
manifold, i.e., the warping function $f$ is not constant on any non empty open subset of $I$.
Moreover, we also assume $\log f$ is convex, which gets that $f$ is also convex. Thus, the
Riemannian manifold $I\times_f M$ admits a global convex function. This kind of functions were widely used, but mostly along curves or locally defined until the significant paper of R.L. Bishop and B. O'Neill \cite{B-O} was published. In this paper, the authors extensively studied how the existence of a convex function affects the topology and the curvature of a Riemannian manifold. Moreover,
convex functions are used to construct Riemannian warped products of negative sectional curvature \cite[Th. 7.5]{B-O}. On the other hand, under suitable assumptions on the fiber of $I\times_f M$,
the condition $(\log f)''\geq 0$ may be interpreted in terms of curvature (see Remark \ref{Remark}(b)).

As in \cite{Al-Da} and \cite{Rosenberg_2002} we will show that a complete hypersurface becomes parabolic under several assumptions. In these references the well-known strong relation between Gauss
curvature and parabolicity in dimension two is used in the minimal case. However, we
will follow here an analogous procedure to obtain parabolicity for a pointwise conformally
related metric to the induced one on a complete hypersurface as in \cite{RRS2} without assuming
before minimality of the hypersurface. Recall that, for higher dimensions there are no
clear relation with assumptions on the sectional curvature and parabolicity (for instance,
$\mathbb{R}^n$ is parabolic if and only if $n\leq 2$). Thus, the fiber $(M^n,g)$ of $I\times_f M^n$ will be assumed to be parabolic. Then, the strategy we will follow along this paper is to show that the parabolicity of the fiber provides the parabolicity of certain Riemannian metric pointwise
conformally related to the induced one on an entire graph in $I\times_f M^n$ , when some natural
assumptions are fulfilled (Lemma \ref{parabolic}). Next, a distinguished positive function on the entire graph is considered and its Laplacian respect to the conformal metric is computed (see
\cite{Osserman}). Under several natural assumptions, this function becomes super-harmonic, leading
that it is constant.

\vspace{1mm}

 First of all, we get (Theorem \ref{propio}),

\vspace{1mm}

\begin{quote}
{\it Let $(M,g)$ be an $n (\geq 2)$-dimensional parabolic Riemannian manifold and $f\in C^{\infty}(I)$, $I \subseteq \R$, positive, non-locally constant and which satisfies $(\log f)''\geq 0$. The only entire solutions $u\in C^{\infty}(M)$ to the {\rm MS} equation $(\ref{edp})$, which satisfy $|D u| \leq c \, f(u)$, for some $c \in \mathbb{R}^{+}$, are the constants.}
\end{quote}

When the warping function is monotone, stronger results are obtained (Subsection \ref{monotone}).
For instance we prove (Theorem \ref{signada2}),

\begin{quote}
{\it Let $(M,g)$ be an $n (\geq 2)$-dimensional parabolic Riemannian manifold and $f\in C^{\infty}(I)$, $I=(a,b) \subseteq \R$, positive and monotone.

\par
\vspace{2mm} \noindent i) If f is non-increasing and
$f\in L^1(a)$, or

\vspace{2mm} \noindent ii) If f is non-decreasing and $f\in L^1(b)$, \vspace{2mm}
\par
\noindent then, the only entire solutions $u\in C^{\infty}(M)$ to the {\rm MS} equation $(\ref{edp})$, whose gradient satisfies $|D u| \leq c \, f(u)$, for some $c \in \mathbb{R}^{+}$, are the
constants.}
\end{quote}

In another setting, under certain boundedness assumption of the Ricci curvature we can
give new results to the minimal hypersurface equation (2) on some parabolic Riemannian
manifolds. (Theorem \ref{tgeodesic})

\begin{quote}
{\it Let $(M,g)$ be an $n (\geq 2)$-dimensional parabolic Riemannian manifold with non-positive definite Ricci tensor. If $u\in C^{\infty}(M)$ is an entire solution to the minimal hypersurface equation $(\ref{edp2.1})$, whose gradient satisfies $|D u| \leq c$, for some $c \in \mathbb{R}^{+}$, then the graph $\Sigma_u$ is totally geodesic in $I\times M$. Moreover, if the Ricci tensor is negative definite at some point $p_0$, then $u$ must be constant.}
\end{quote}

It should be noted that, in the very particular but important case  $ M^2=\mathbb{R}^2,$, this
result provides a new proof of the Bernstein theorem (Remark \ref{bers}). Finally, when $M^2$ is
a complete cylinder on an $n(\geq 2)$-dimensional compact Riemannian manifold $M^n,g)$, the
topology of $M^n$ may be used to characterize all the solutions of the minimal hypersurface
equation (\ref{edp2.1}). (Theorem \ref{moser}),

\begin{quote}
{\it Let $(M^n,g)$ be an $n$-dimensional compact Riemannian manifold with non-positive definite Ricci tensor and assume the Euler-Poincar\'e characteristic of $M$ is non zero. The only entire
solutions $u$ to the minimal hypersurface  equation $(\ref{edp2.1})$ on $\R \times M$,
with bounded length of its gradient, are the functions $u(s,x)=as+b$, $a,b\in \R$.}
\end{quote}

\section{Preliminaries}
Let us consider the warped product $\overline{M}:=I \times_f M$ of base $(I,dt^2)$, fiber the Riemannian manifold $(M,g)$ and warping function $f$. The vector field $K:=f(\pi_{_I}) \,\partial_t$,
where $\partial_t$ is the coordinate vector field, satisfies \cite[Prop.
7.35]{ON},
\begin{equation}\label{basica}
\overline{\nabla}_X K=f'(\pi_{_I}) \, X,
\end{equation}
for any $X\in\mathfrak{X}(\overline{M})$, where $\overline{\nabla}$ is the Levi-Civita connection of $\bar g$. Thus, the vector field $K$ is conformal, with $\mathcal{L}_K{\bar g}=2\,f'(\pi_{_I})\,{\bar g}$, its metrically equivalent 1-form is closed and  its divergence satisfies
$\mathrm{div}(K)=(n+1)f'(\pi_{_I})$.

\vspace{1mm}

For each $u\in C^{\infty}(M)$ let $\Sigma_u=\lbrace (u(p),p)\, : \, p\in M\rbrace$ be the entire graph defined by $u$ on $M$. The subset $\Sigma_u$ is a regular hypersurface in $\overline{M}$ and
it inherits a Riemannian metric $g_{_{\Sigma_u}}$ from $\overline{M}$ which, on $M$, has the following expression,
\begin{equation} \label{metrica}
g_u = du^2 + f(u)^2 \, g \, ,
\end{equation}

\noindent where $f(u):=f\circ u$. If we put $\tau:=\pi_{_I}\circ i$, where $i$ is the inclusion of
$\Sigma_u$ in $\overline{M}$, then it is no difficult to obtain that the gradient of $\tau$ satisfies
\begin{equation}\label{gradiente}
\nabla \tau=\,{\partial_t}^{T},
\end{equation}
where ${\partial_t}^T$ denotes the projection of $\partial_t$ on $\Sigma_u$.

\vspace{1mm}

A unit normal vector field of $\Sigma_u$ in $\overline{M}$ is
\begin{equation}\label{unit}
N=\frac{f(u)}{\sqrt{f(u)^2+\vert
Du\vert^2}}\,\Big(\partial_t-\frac{1}{f(u)^2}\,Du\Big),
\end{equation}
where $Du$ is the gradient of the function $u$ in $(M,g)$, and $\vert Du\vert^2:=g(Du,Du)$.
Clearly,
\begin{equation} \label{coseno}
\cos \theta=\frac{f(u)}{\sqrt{f(u)^2+\vert Du\vert^2}}\, ,
\end{equation}
where $\theta$ is the angle between $N$ and $\partial_t$.
Therefore, from (\ref{gradiente}) we have
\begin{equation}\label{nueva}
\|\nabla\tau\|^2:=g_{_{\Sigma_u}}(\nabla\tau,\nabla \tau)=\sin^2\theta.
\end{equation}

\vspace{1mm}

The Gauss and Weingarten formulas of $\Sigma_u$ in $\overline{M}$ are respectively written
\begin{equation}
\overline{\nabla}_X Y=\nabla_X Y+g_{_{\Sigma_u}}(AX,Y)\, N
\end{equation}
\begin{equation}
AX=-\overline{\nabla}_X N
\end{equation}
for all $X,Y\in\mathfrak{X}(\Sigma_u)$, where $\nabla$ is the Levi-Civita connection of the induced metric on $\Sigma_u$ and $A$ is the shape operator associated to $N$. The mean curvature
function relative to $N$ is $H=\frac{1}{n}\mathrm{trace}(A)$.  As it is well-known, $H=0$ if and only if $\Sigma_u$ is locally a critical point of the $n$-dimensional volume functional for
compactly supported normal variations. The graph $\Sigma_u$ is said to be minimal when $H=0$.

In $\overline{M}$ the graphs of any constant function $u=t_0$, $t_0\in I$ (i.e., the level hypersurfaces of the projection $\pi_{_I}\longrightarrow I$) constitute a distinguished family of hypersurfaces in $\overline{M}$, the so-called slices $t=t_0$. The normal vector field $N$ of a slice $t=t_0$ is the restriction of $\partial_t$ to $t=t_0$. From (\ref{basica}), the shape operator with respect to $N$ is given by $A=- f'(t_0)/f(t_0) I$, where $I$ denotes the identity transformation. Therefore, a slice $t=t_0$ is a totally umbilical hypersurface of constant mean curvature
\begin{equation}\label{meancurvature}
H = - f'(t_0)/f(t_0).
\end{equation}

Thus, a slice $t=t_0$ is minimal if and only if $f'(t_0)=0$ (i.e., if and only if it is totally geodesic). Note that a slice $t=t_0$, with $f'(t_0)=0$ gives a trivial entire solution to the minimal hypersurface equation (\ref{edp2.1}). Under several geometric assumptions we will prove that these
slices provide the only entire solutions to minimal hypersurface equation (\ref{edp2.1}).

\vspace{1mm}

Coming back to an arbitrary graph $\Sigma_u$ in $\overline{M}$, consider the tangential component
$K^T:=K-\bar{g}(N,K)N$ on $\Sigma_u$ of $K$. From (\ref{basica}) and using the Gauss and Weingarten formulas we get

$$
\nabla_XK^T=f'(\tau)\, X+f(\tau)\, {\bar g}(N,\partial_t)\, AX
$$

\noindent for any $X\in \mathfrak{X}(\Sigma_u)$. Making use of (\ref{gradiente}), the Laplacian of $\tau$ on $\Sigma_u$ is given by
\begin{equation}\label{gradient1}
\Delta\tau=\frac{f'(\tau)}{f(\tau)}\lbrace
n-\|\nabla\tau\|^2\rbrace +nH \, \bar g(N,\partial_t).
\end{equation}
Observe that $n- \|\nabla \tau\|^2>0$ by (\ref{nueva}) and $n\geq 2$. A direct computation from (\ref{gradiente}) and (\ref{gradient1}) gives

\begin{equation}\label{laplacian2}
\Delta f(\tau) = n \,\frac{f'(\tau)^2}{f(\tau)} +f(\tau) ( \log
f)''(\tau) \|\nabla \tau\|^2+nH \, f'(\tau) \, \bar g(N,\partial_t).
\end{equation}

\section{The minimal hypersurface equation}
From Weingarten formula and taking into account \cite[Prop. 7.35]{ON}, the shape operator $A$ of a graph $\Sigma_u=\{(u(p),p)\, :\, p\in M^n\}$ of $I\times _f M^n$, corresponding to $N$ given in (\ref{unit}), satisfies,

$$A(X)=\frac{-f'(u)}{\sqrt{f(u)^2+\mid Du\mid^2}}\Big\{\frac{f'(u)}{f(u)}X+\frac{f'(u)\, g(Du,X)}{f(u)\sqrt{f(u)^2+\mid Du\mid^2}}Du$$

$$-\frac{g(D_X Du,Du)}{f(u)^2(f(u)^2+\mid Du\mid^2)}Du-\frac{1}{f(u)^2}D_XDu\Big\},$$

\noindent for all $X$ tangent to the graph. The contraction of this formula when $H=0$ leads to the minimal hypersurface equation (\ref{edp2.1}).

Alternatively, let $\Omega$ be an open domain in a Riemannian manifold $M^n$, let $f:I\longrightarrow\mathbb{R}$ be a positive smooth function and let $u: \Omega
\rightarrow \mathbb{R}$ be a smooth function such that $u(\Omega)\subset I$. The volume of the graph restricted to a compact subset $Q$ in $\Omega$ is computed as follows,

\begin{equation}
\mathrm{vol}(\Sigma_u, Q) = \int_Q f(u)^{n-1}\, \sqrt{f(u)^2 + |D
u|^2} \; d\mu_{g_u} \ ,
\end{equation}

\noindent where $d\mu_{g_u}$ is the canonical measure associated to $g_u$

Consider a smooth function $v: \Omega \rightarrow \mathbb{R}$ with compact support $Q$ in $\Omega$. The volume of the graph of the function $u + tv$, $t \in \R$, is
\begin{equation} \displaystyle
\int_Q f(u +tv)^{n-1} \sqrt{f(u+tv)^2 + |Du+tDv)|^2} \; d\mu_{g_u}
\, .
\end{equation}
Assume
\begin{equation}\label{derivadavolumen}
\frac{d}{dt}\Bigl|_{t=0} \mathrm{vol}(\Sigma_{u+tv},Q)=0
\end{equation}
for every smooth function $v$ on $\Omega$ with compact support. A standard argument from (\ref{derivadavolumen}) gets (\ref{edp2.1}).

\section{Parabolicity in Riemannian manifolds}
A non-compact Riemannian manifold is said to be parabolic if it admits no non-constant positive superharmonic function (see \cite{Kazdan}, for instance). In the two dimensional case, this
notion is very close to the classical parabolicity for Riemann surfaces. Moreover, it is strongly related to the behaviour of the Gauss curvature of the surface. First of all, the seminal result
by Ahlfors and Blanc-Fiala-Huber states that a complete 2-dimensional Riemannian manifold with non-negative Gauss curvature is parabolic (see \cite{H}, \cite{Kazdan}). There are another results in this direction, for example if the Gauss curvature of a complete Riemannian surface satisfies $K\geq -1/(r^2\log r)$, for $r$, the distance to a fixed point, sufficiently large, then the surface is parabolic \cite{GW}. And if a complete Riemannian surface is such that the negative part of its Gauss curvature is integrable, then the surface must be
parabolic \cite{Li}.

\vspace{1mm}

For higher dimensions, parabolicity of Riemannian manifolds has a di different behaviour and, in particular, it has no clear relation with assumptions on the sectional curvature. In fact, the
Euclidean space $\R^n$ is parabolic if and only if $n\leq 2$. Even more, if $(M_1,g_1)$ is any compact Riemannian manifold and $(M_2,g_2)$ is a parabolic Riemannian manifold, then $M_1\times
M_2$ endowed with the product metric $g_1 + g_2$ is parabolic, \cite{Kazdan}. In particular, the product of a compact Riemannian manifold and the real line $\R$ is always parabolic. On the other hand, parabolicity is closely related with the volume growth of the geodesic balls in an $n(\geq 2)$-dimensional non-compact complete Riemannian manifold $(M,g)$; indeed, if it has moderate volume growth, then $(M,g)$ must be parabolic \cite{Kar}.

An important property of parabolicity is that it is invariant under quasi-isometries \cite[Cor. 5.3]{GR}, \cite{Ka2}. Let us recall that given Riemannian manifolds $(P,g)$ and $(P',g')$, a
diffeomorphism $\varphi$ from $P$ onto $P'$ is called a quasi-isometry if there exists a constant $c \geq 1$ such that 
$$c^{-1}\vert v\vert_{g}\,\leq\,\vert d\varphi(v)\vert_{g'}\,\leq\, c
\, \vert v\vert_{g}$$ \
\noindent for all $v\in T_p P$, $p\in P$ (see \cite{Ka} for more details). Moreover, this result can be used to construct new parabolic Riemannian manifolds as follows. Consider a parabolic Riemannian manifold $(M,g)$ and let $h \in C^\infty(M)$ such that $\inf(h)>0$ and $\sup(h)<\infty$. Then, the Riemannian manifold $(M, h^2 \, g)$ is quasi-isometric to $(M,g)$ and, therefore, it is also parabolic. On the other hand, suppose that $(M_1, g_1)$ and $(M_2, g_2)$ are parabolic Riemannian manifolds such that $(M_1 \times M_2, g_1 + g_2)$ is also parabolic (of course, a Riemannian product of parabolic manifolds is not parabolic in general). For any $h \in C^\infty(M_1)$ such that $\inf (h)>0$ and $\sup(h)<\infty$, we have that $(M_1 \times M_2, g_1 + h^2 \, g_2)$ is parabolic. In fact,
writing $c = \inf(h)$ and $d = \sup(h)$, the following inequalities holds,
\begin{eqnarray*}
(g_1 + h^2 \, g_2)(X,X) &\leq & g_1 (X_1,X_1) + d^2 \, g_2 (X_2, X_2) \\
 & \leq & (1+d^2) (g_1 + g_2) (X,X) \, .
\end{eqnarray*}
\begin{eqnarray*}
(g_1 + h^2 \, g_2)(X,X) & \geq & g_1 (X_1,X_1) + c^2 \, g_2 (X_2,X_2) \\
                    & \geq & \min \{ 1,c^2 \} \, (g_1 + g_2) (X,X) \,
                    ,
\end{eqnarray*}
where $X=(X_1,X_2)$, which mean that $(M_1\times M_2,g_1+g_2)$ and $(M_1\times M_2,g_1+h^2 \,g_2)$ are quasi-isometric. Observe that the same argument shows that if $(M_1, g_1)$ is a compact
Riemannian manifold, $(M_2, g_2)$ a parabolic Riemannian manifold and $h \in C^\infty(M_1)$, $h>0$, then\linebreak $(M_1\times M_2,g_1+h^2\,g_2)$ is also parabolic.

\section{Main results}
We begin this section with the statement of a technical result to get the parabolicity of certain conformal metric  of an entire graph $\Sigma_u$ in $\overline{M} = I \times_f M$ from the parabolicity of the fiber.

\begin{lema}\label{parabolic} Let $(M,g)$ be a parabolic Riemannian manifold and let $f$ be a positive smooth function on the interval $I$. If $u \in C^\infty (M)$ satisfies $u(M^n)\subset I$ and $|Du | \leq C \, f(u)$, for some $c \in \R^{+}$, for some $c\in \mathbb{R}^+$, then the metric
\begin{equation}\label{conmetric01}
\widehat{g}:= \frac{1}{f(u)^2}\, g_u
\end{equation}
\noindent where $g_u$ is given in $(\ref{metrica})$, is also parabolic on $M^n$; i.e., the graph $\Sigma_u$ endowed with the metric
\begin{equation}
\widetilde{g}:=\frac{1}{f(\tau)^2}\,g_{_{\Sigma_u}}
\end{equation}\label{conmetric}
\noindent is parabolic.
\end{lema}
\noindent\emph{Proof.} From (\ref{metrica}) we easily get

\begin{equation} \label{des1}
\widehat{g}(X,X) \geq g(X,X) \, ,
\end{equation}
for any $X \in \mathfrak{X}(M)$. Now, taking into account the Schwarz inequality for $g$, we have
$$
\widehat{g} (X,X) \leq \left( 1 + \frac{|D u|^2}{f(u)^2} \right) \, g(X,X) \, ,
$$ which turns into

\begin{equation} \label{des2}
\widehat{g}(X,X) \leq (1 + C^2 ) \,  g(X,X) \, ,
\end{equation}
from our assumption. Therefore, from (\ref{des1}) and (\ref{des2}) we conclude that the identity map is a quasi-isometry from $(M,\widehat{g})$ onto $(M,g)$ which ends the proof. \hfill{$\Box$}

\begin{rema} {\rm The assumption on $\mid Du\mid$ in lemma \ref{parabolic} has a clear geometrical meaning. In fact, from (\ref{coseno}), the angle between the coordinate vector field $\partial_t$ and the unit normal vector field $N$, given in (\ref{unit}), is bounded away from $\pi/2$.}
\end{rema}

Using again the invariance by quasi-isometries of parabolicity, we get,

\begin{coro}\label{consecuencia1}
Let $(M^n,g)$ be a parabolic Riemannian manifold and let $f$ be a
positive smooth function on $I$ such that $\inf f >0$. If $u \in C^\infty (M)$ satisfies $u(M^n)\subset I$ and $|Du | \leq C \, f(u)$, for some $c \in \R^{+}$, for some $c\in \mathbb{R}^+$,
then the metric $g_u$, given in (\ref{metrica}), is parabolic on $M^n$.
\end{coro}

\begin{coro}\label{consecuencia2}
Let $(M^n,g)$ be a parabolic Riemannian manifold. If $u \in C^\infty
(M)$ satisfies $| Du | \leq c$, for some $c \in \R^{+}$, then the product metric
$g_u= du^2+g$ is parabolic on $M^n$.
\end{coro}

Now observe that the equations (\ref{gradient1}) and (\ref{laplacian2}) may be rewritten for a minimal hypersurface in term of the conformal metric $\tilde{g}$, given in (\ref{conmetric}), as follow,

\begin{equation} \label{conforme1}
\widetilde{\Delta} \tau = n\,f(\tau) \, f'(\tau) - (n-1)\, \frac{f'(\tau)}{f(\tau)}\,||{\widetilde \nabla}
\tau||_{\widetilde g}^2 ,
\end{equation}

\begin{equation} \label{conforme2}
\widetilde{\Delta} f(\tau) = n f(\tau) \, f'(\tau)^2 + \left\{ f(\tau)
\, (\log f)''(\tau) -(n-2)\, \frac{f'(\tau)^2}{f(\tau)} \right\} \, ||{\widetilde \nabla}
\tau||_{\widetilde g}^2.
\end{equation}

\vspace{1mm}

Now, we can give some uniqueness results using previous work in this section.

\begin{teor} \label{propio}
Let $(M^n,g)$ be an $n (\geq 2)$-dimensional parabolic Riemannian manifold and $f\in C^{\infty}(I)$, $I \subseteq \R$, positive, non-locally constant and which satisfies $(\log f)''\geq 0$. The only entire solutions $u\in C^{\infty}(M)$ to the {\rm MS} equation $(\ref{edp})$, which satisfy $|D u| \leq c \, f(u)$, for some $c \in \mathbb{R}^{+}$, are the constants.
\end{teor}

\noindent\emph{Proof.} First of all, Lemma \ref{parabolic} may be called to obtain that $M^n$ endowed with the conformal metric $\widehat{g}=(1/f(u)^2)g_u$ given in (\ref{conmetric01}) is parabolic.  Now, using equation (\ref{conforme2}), we have that the $\widetilde{g}$-Laplacian of the positive and bounded function
$$\mathrm{arccot}\, f(\tau):M \rightarrow [0,2\pi)$$
\noindent satisfies,

\begin{eqnarray*}
\widetilde{\Delta}\, \mathrm{arccot} f (\tau) & = &
-\frac{f(\tau)}{1+f(\tau)^2}\,(\log f)''(\tau)\,||{\widetilde
\nabla}
\tau||_{\widetilde g}^2 \\[2mm]
{}&&-
\frac{f'(\tau)^2}{f(\tau)^2\left(1+f(\tau)^2\right)}\Big\{n(f(\tau)^2-
||{\widetilde \nabla} \tau||_{\widetilde g}^2)f(\tau)
 + \frac{2}{f(\tau)\left(1+f(\tau)^2\right)}\,||{\widetilde \nabla}
\tau||_{\widetilde g}^2\, \Big\} \, .
\end{eqnarray*}

\noindent Taking into account (\ref{nueva}) we get $||{\widetilde \nabla} \tau||_{\widetilde g}^2 < f(\tau)^2$. From the assumptions, we have that the function ${\rm arccot} f(\tau)$ is $\widetilde{g}$-superharmonic. Then, as a consequence of the $\widetilde{g}$-parabolicity of $\Sigma_u$,  $f(\tau)$ must be constant and, consequently, $\tau$ is constant.\hfill{$\Box$}

\begin{rema}\label{Remark}{\rm {\bf (a)} The assumption ``non-locally constant'' on $f$ cannot be removed clearly as the entire solutions to the MS equation (\ref{edp2}) shows. {\bf (b)} The assumption $(\log f)'' \geq 0$ has also a clear geometrical meaning. It guarantees that
the Ricci curvature of $I\times_f M$ in non-positive whenever the Ricci curvature of the fiber is non-positive \cite[Cor. 7.43]{ON}. On the other hand, taking into account (\ref{meancurvature}), it
implies that the mean curvature of the leaves is a non-increasing function.}
\end{rema}

\begin{coro}\label{propio2}
Let $(M^n,g)$ be an $n (\geq 2)$-dimensional parabolic Riemannian manifold and $f\in C^{\infty}(I)$, $I \subseteq \R$, non-locally constant and such that $\inf f >0$  and $(\log f)''\geq 0$. The only entire solutions $u\in C^{\infty}(M)$ to the {\rm MS} equation $(\ref{edp})$, which satisfy $|D u| \leq c$, for some $c \in \mathbb{R}^{+}$, are the constants.
\end{coro}

For the case of a more general warping function, we have,

\begin{teor} \label{nopropio}
Let $(M^n,g)$ be an $n (\geq 2)$-dimensional parabolic Riemannian manifold and $f\in C^{\infty}(I)$, $I \subseteq \R$, positive and such that $(\log f)''\geq 0$. The only entire solutions $u\in
C^{\infty}(M)$ to the {\rm MS} equation $(\ref{edp})$, which satisfy $|Du| \leq c \, f(u)$, for some $c \in \mathbb{R}^{+}$, and which are bounded from below or from above, are the constants.
\end{teor}

\noindent\emph{Proof.} We already know by the previous theorem that $f(\tau)$ is constant, now (\ref{laplacian2}) implies $f'(\tau)=0$. Therefore, from (\ref{gradient1}), $\tau$ is harmonic on $M$,
ending the proof. \hfill{$\Box$}

\vspace{1mm}

In the the case of $f=1$, i.e., for an ambient Riemannian product, the previous theorem specializes to the following result,

\begin{coro}\label{nopropioconsecuencia}
Let $(M,g)$ be an $n (\geq 2)$-dimensional parabolic Riemannian manifold. The only entire solutions $u\in C^{\infty}(M)$ to the {\rm MS} equation $(\ref{edp2.1})$, which satisfy $|D u| \leq C$, for some $C
\in \mathbb{R}^{+}$, and which are bounded from below or from above, are the constants.
\end{coro}

\begin{counter}\label{contra1} {\rm \textbf{(a)} Consider the family of complete Riemannian manifolds $\mathbb{R}\times_f\mathbb{R}$ (completeness follows from \cite[Lem. 7.40]{ON}). A function $u=u(x)$
is solution of the MS equation (\ref{edp2}) on $\mathbb{R}\times_f\mathbb{R}$ if and only if it satisfies
$$
f(x)\, \frac{u'(x)}{\sqrt{1+u'(x)^2}}= C \,
$$
for some $C \in \mathbb{R}$ and any $x \in \mathbb{R}$. If we choose $f(x)=\sqrt{1 + \cosh^4 (x)}$, then the corresponding metric $dx^2+f(x)^2dy^2$ is not parabolic (the function $v(x)=-1/\cosh^2x$, $x\in \mathbb{R}$, satisfies $\Delta v=2\left[(\cosh^2x-1)^2+2\right]/\left[\cosh^4x(1+\cosh^4x)\right]<0$).
Now, the function $$u(x) = \tanh x $$ is a solution of the MS equation (\ref{edp2}) on $\mathbb{R}\times_f\mathbb{R}$. It is trivially bounded and the norm of its gradient is found also to be bounded. \textbf{(b)} Now consider the complete Riemannian manifold $\mathbb{R}\times_h\mathbb{R}$ where $h(x)=(\sqrt{2x^4+6x^2+5})/(x^2+2)$. As the warping function $h$ satisfies $\sqrt{5}/2 \leq h(x) < \sqrt{2}$, at any $x\in \mathbb{R}$, from the considerations showed in Section $4$ we
conclude that $\mathbb{R}\times_h\mathbb{R}$ is parabolic. The function $$w(x) = x+\arctan x$$ is a solution of the MS equation (\ref{edp2}) on $\mathbb{R}\times_h\mathbb{R}$, with bounded
length of its gradient. Note that $w$ is unbounded neither from below nor from above.}
\end{counter}

\subsection{The case of monotone warping function}\label{monotone}
Recall that a positive continuous function $f$ on $(a,b)$, $-\infty \leq a<b\leq \infty$, is said to satisfy $f\in L^1(a)$ (resp. $f\in L^1(b)$) if $\int_a^c f(s)ds<\infty$ (resp. $\int_c^b
f(s)ds<\infty$) for some $c\in(a,b)$.
\begin{teor} \label{signada1}
Let $(M^n,g)$ be an $n (\geq 2)$-dimensional parabolic Riemannian manifold and $f\in C^{\infty}(I)$, $I=(a,b) \subseteq \R$,
positive and monotone.

\par
\vspace{2mm} \noindent i) If $a \in  \mathbb{R}$ and f is non-increasing with
$f\in L^1(a)$, or

\vspace{2mm} \noindent ii) If $b \in  \mathbb{R}$ and f is non-decreasing with 
$f\in L^1(b)$, \vspace{2mm}
\par
\noindent then, the only entire solutions $u\in C^{\infty}(M)$ to the {\rm MS} equation $(\ref{edp})$, whose gradient satisfies $|D u| \leq c \, f(u)$, for some $c \in \mathbb{R}^{+}$, are the
constants.
\end{teor}

\noindent\emph{Proof.} Under the assumption {\it i)}, we consider the function $\mathcal{F}(\tau)$ on $\Sigma_u$, defined by
$$
\mathcal{F}(\tau) = \int_{s_0}^{\tau} f(s) \, ds \, ,
$$where $s_0 = \inf \left(\tau\right)$. Clearly, $\mathcal{F}(\tau) \geq 0$ and the Laplacian of $\mathcal{F}(\tau)$ respect to the conformal metric $\widetilde{g}$ satisfies
\begin{equation}
\widetilde{\Delta} \mathcal{F}(\tau) =  f'(\tau) f(\tau)^2
\left\{n -(n-2)\,\frac{||\widetilde{\nabla} \tau||_{\widetilde
{g}}^2}{f(\tau)^2} \right\}\,\leq 0\, \, . \label{laplacianoF}
\end{equation}
Using now the parabolicity of $\widetilde{g}$, we get that $u$ must be constant.

\vspace{1mm}

The case {\it ii)} follows analogously changing
$\mathcal{F}(\tau)$ to the function
$\widetilde{\mathcal{F}}(\tau)$ given by
$$ \widetilde{\mathcal{F}}(\tau) = \int_{\tau}^{s^0} f(s) \, ds \,
,
$$where $s^0 = \sup (\tau)$.\hfill{$\Box$}

\begin{coro} \label{signada1consecuencia1}
Let $(M^n,g)$ be an $n (\geq 2)$-dimensional parabolic Riemannian manifold and $f\in C^{\infty}(I)$, $I=(a,b) \subseteq \R$, positive, monotone and $\inf{f} > 0$.
\par
\vspace{2mm}
\noindent i) If $a \in  \mathbb{R}$ and f is non-increasing with
$f\in L^1(a)$, or

\vspace{2mm} \noindent ii) If $b \in  \mathbb{R}$ and f is non-decreasing with 
$f\in L^1(b)$,
\vspace{2mm}
\par
\noindent then, the only entire solutions $u\in C^{\infty}(M)$ to the {\rm MS} equation $(\ref{edp})$, whose gradient satisfies $|D u| \leq c$, for some $c \in \mathbb{R}^{+}$, are the constants.
\end{coro}

Observe that if only bounded from below (or from above) solutions to the MS equation (\ref{edp}) are considered, the $L^1$ assumptions on $f$ can be dropped.

\begin{teor} \label{signada2}
Let $(M,g)$ an $n (\geq 2)$-dimensional parabolic Riemannian manifold and $f\in C^{\infty}(I)$, $I \subseteq \R$, positive and non-increasing $($resp. non-decreasing$)$. The only entire solutions
$u\in C^{\infty}(M)$ to the {\rm MS} equation $(\ref{edp})$, bounded from below $($resp. bounded from above$)$, whose gradient satisfies $|Du| \leq C \, f(u)$, for some $C \in \mathbb{R}^{+}$, are the constants.
\end{teor}

We can give a nice direct consequence of the last theorem,

\begin{coro}\label{signada2consecuencia}
Let $(M,g)$ be an $n (\geq 2)$-dimensional parabolic Riemannian manifold and $f\in C^{\infty}(I)$, $I \subseteq \R$, positive and monotone. The only entire bounded solutions $u\in C^{\infty}(M)$ to the {\rm MS} equation $(\ref{edp})$, whose gradient satisfies $|D u| \leq C$, for some $C \in \mathbb{R}^{+}$, are the constants.
\end{coro}

Next, in order to make use of a basic result on the $\varphi$-Laplacian, \cite[Ch. 5]{pigola}, we will transform MS equation (\ref{edp}) by means of a suitable change of variable. In fact, let $u \in C^\infty(M)$ an entire solution to the MS equation and define
$$
v= \psi(u),\quad \text{where} \quad \psi(t):=\int_{u_0}^{t}
\frac{1}{f(s)} \, ds \, ,
$$ and $u_0$ is some value of $u$. Taking into account $f(u) \, D v = D u$, equation (\ref{edp}) can be written, in terms of $v$ as follows,
\begin{equation}\label{transformada}
\mathrm{div} \left( \frac{D v}{\sqrt{1+ |D v|^2}} \right) = n
\,\frac{ f'((\psi^{-1})(v))}{\sqrt{1+|Dv|^2}}  \, ,
\end{equation}
whose right hand side can be seen as the $\varphi$-Laplacian of $v$, where $\varphi \in C^\infty([0,\infty))$ is given by $\varphi (x) = x/\sqrt{1+x^2}$. Now, we can state,

\begin{teor}
Let $(M,g)$ be an $n(\geq 2)$-dimensional complete Riemannian manifold with quadratic volume growth, and let $f\in C^{\infty}(I)$, $I \subseteq \R$, a smooth positive non-increasing
$($resp. non-decreasing$)$ function. The only entire bounded below $($resp. above$)$ solutions to the {\rm MS} equation $(\ref{edp})$, are the constants.
\end{teor}

\noindent\emph{Proof.} Let $u \in C^\infty(M)$ an entire solution to the MS equation (\ref{edp}). Up reversing the $t$ coordinate if it is necessary, we only consider the case that $u$ is bounded below and $f'\leq 0$. Now, making use of (\ref{transformada}), the transformed function $v$
of $u$ is a bounded below $\varphi$-subharmonic function on the complete Riemannian manifold $(M,g)$ which has quadratic volume growth. The result follows directly from \cite[Th. 5.1]{pigola}.
\hfill{$\Box$}

\vspace{2mm}

As a particular case,
\begin{coro}
Let $(M,g)$ be an $n(\geq 2)$-dimensional complete Riemannian manifold with quadratic volume growth. The only entire bounded above or below solutions to the {\rm MS} equation $(\ref{edp2})$, are the constants.
\end{coro}

Observe that the integral assumption needed in \cite[Th. 5.1]{pigola} is satisfied for a wide family of Riemannian manifolds bigger than the quadratic volume growth ones. Moreover, for the considered $\varphi$-Laplacian, this assumption is in fact equivalent to condition (7.15) in \cite[Th. 7.5]{GR}. In the same reference, it is shown that this is also a necessary condition for the
parabolicity for spherically symmetric manifolds. Therefore, we can derive the following consequence,

\begin{coro}
Let $(M,g)$ be an $n(\geq 2)$-dimensional spherically symmetric parabolic Riemannian manifold, and let $f\in C^{\infty}(I)$, $I \subseteq \R$, a smooth positive non-increasing $($resp. non-decreasing$)$
function. The only entire bounded below $($resp. above$)$ solutions to the {\rm MS} equation $(\ref{edp})$, are the constants.
\end{coro}

\subsection{The case of constant warping function}

When a product ambient space $I\times M$ is considered, some extra curvature assumption on the Riemannian manifold $(M,g)$, besides its parabolicity, is needed to arrive to a result of uniqueness.

\begin{teor}\label{tgeodesic}
Let $(M,g)$ be an $n (\geq 2)$-dimensional parabolic Riemannian manifold with non-positive definite Ricci tensor. If $u\in C^{\infty}(M)$ is an entire solution to the {\rm MS} equation $(\ref{edp2})$, whose gradient satisfies $|D u| \leq C$, for some $C \in \mathbb{R}^{+}$, then the graph $\Sigma_u$ is totally geodesic in $I\times M$. Moreover, if the Ricci tensor of $M$ is negative definite at some point $p_0$, then $u$ must be constant.
\end{teor}

\noindent\emph{Proof.} Consider the function $\cos \theta$ on the graph $\Sigma_u$ and note that
$$
\nabla \cos \theta = - A \, \nabla \tau \, ,
$$from (\ref{basica}) and (\ref{gradiente}), and consequently
$$
\Delta \cos \theta=-\,\mathrm{div}_{_{\Sigma_{_u}}}(A\nabla\tau).
$$Now the Gauss and Codazzi equations can be claimed to get,
\begin{equation}\label{lcos}
\Delta\cos\theta=-\cos \theta \, \mathrm{Ric}^M(N^M,N^M)+\cos\theta\,
\mathrm{trace}(A^2),
\end{equation}where $N^M$ means the projection onto $M$ of the unit normal vector field $N$. Now, the use of the curvature assumption in (\ref{lcos}) gives that $\cos \theta$ is subharmonic. Using now
the parabolicity of $(M,g)$, we conclude that $\theta$ must be constant. Consequently, $A=0$. Finally, the last assertion follows taking into account that if the function $\cos\theta$ is constant,
then $\cos \theta(p_0)=0$, which implies that $u$ must be also constant. \hfill{$\Box$}

\vspace{1mm}

\begin{rema}\label{bers}{\rm Coming back to the classical Bernstein theorem for $n=2$, observe
that an entire solution $u$ to the MS equation (\ref{edp2}) satisfies $|Du| \leq C$, for some $C\in \R^+$ according Bers' theorem \cite{Bers}. Therefore, $u$ lies under the assumptions of
Theorem \ref{tgeodesic} and consequently $\Sigma_u$ must be a plane in Euclidean space $\R^3$. Alternatively, it is possible to give another argument, using also Theorem \ref{tgeodesic}, leading to the classical Bernstein theorem. Namely, let $\Sigma_u$ be the graph of an entire solution $u$ to the MS equation (\ref{edp2}). The parabolicity of the graph can be deduced directly as follows. Let
$p_0\in\Sigma_u$ be and consider the closed ball $\bar B_{\R^3}(p_0,r)\subset\R^3$, centered at $p_0$, with radius $r$. Clearly
$$\mathrm{Area}\Big(\bar B_{\R^3}(p_0,r)\cap \Sigma_u\Big)\leq 2\pi r^2,$$
because $\Sigma_u$ is a minimizing surface area.
On the other hand,  the geodesic disk $D_{\Sigma_u}(p_0,r)$ in $\Sigma_u$ satisfies
$$D_{\Sigma_u}(p_0,r)\subset \bar B_{\R^3}(p_0,r)\cap \Sigma_u,$$

\noindent and as a consequence

$$\mathrm{Area}(D_{\Sigma_u}(p_0,r))\leq 2\pi r^2,$$
which in particular means that the area of the geodesic disk has quadratic growth, and this implies that $\Sigma_u$ is parabolic \cite{C-Y}. Now, Theorem \ref{tgeodesic} can be claimed again to get that
$\theta$ is constant, which ends the argument.}
\end{rema}

\vspace{1mm}

\begin{teor}\label{moser}
Let $(F,g_{F})$ be an $(n-1)$-dimensional compact Riemannian manifold with non-positive definite Ricci tensor and assume the Euler-Poincar\'e characteristic of $F$ is non zero. The only entire
solutions $u$ to the {\rm MS} equation $(\ref{edp2})$ on $\R \times F$, with bounded length of its gradient, are the functions $u(s,x)=as+b$, $a,b\in \R$.
\end{teor}

\noindent\emph{Proof.} The Riemannian manifold $M= \R \times F$ is parabolic because $F$ is compact. Corollary \ref{consecuencia2} can be then claimed to obtain that the graph $\Sigma_u$ is parabolic. On the other hand, it is easy to see that the Ricci tensor of $M$ is also non-positive definite. Therefore, as a consequence of Theorem \ref{tgeodesic}, $\Sigma_u$ must be totally geodesic.

The projection $N^M$ of the unitary normal vector field $N$ onto $\R\times F$  is parallel. Therefore, its $F$-component is also parallel on $F$ and vanishes because the Euler-Poincar\'e characteristic of $F$ is non zero. Hence, this leads to $N^M=a \,\partial_s$, for some $a \in \R$. \hfill{$\Box$}

\vspace{1mm}

We end the paper showing that the class of Riemannian manifolds under the assumptions of Theorem \ref{moser} is very wide.

\begin{rema}{\rm
In the previous result, we may take $F$ as the Riemannian product $T^k \times R$, where $T^k$ is a $k$-dimensional flat torus and $R$ is either a compact Ricci-flat Riemannian manifold \cite{F-W}
or a compact negatively Ricci curved Riemannian manifold \cite{G-Y}, \cite{JL}.}
\end{rema}

\end{document}